\documentclass[12pt]{article}
\usepackage{amssymb}
\usepackage{latexsym}
\newcommand{\bit}{\begin{itemize}}
\newcommand{\eit}{\end{itemize}}
\newcommand{\ben}{\begin{enumerate}}
\newcommand{\een}{\end{enumerate}}
\newcommand{\bre}{\begin{rem}}
\newcommand{\ere}{\end{rem}}
\newcommand{\ble}{\begin{lem}}
\newcommand{\ele}{\end{lem}}
\newcommand{\bth}{\begin{thm}}
\renewcommand{\eth}{\end{thm}}
\newcommand{\bmth}{\begin{mainthm}}
\newcommand{\emth}{\end{mainthm}}
\newcommand{\bpr}{\begin{prop}}
\newcommand{\epr}{\end{prop}}
\newcommand{\bco}{\begin{cor}}
\newcommand{\eco}{\end{cor}}

\newcommand{\bcon}{\begin{conj}}
\newcommand{\econ}{\end{conj}}
\newcommand{\bde}{\begin{defn}}
\newcommand{\ede}{\end{defn}}
\newcommand{\bex}{\begin{exa}}
\newcommand{\eex}{\end{exa}}
\newcommand{\barr}{\begin{array}}
\newcommand{\earr}{\end{array}}
\newcommand{\btab}{\begin{tabular}}
\newcommand{\etab}{\end{tabular}}
\newcommand{\beq}{\begin{equation}}
\newcommand{\eeq}{\end{equation}}
\newcommand{\bea}{\begin{eqnarray*}}
\newcommand{\eea}{\end{eqnarray*}}
\newcommand{\bce}{\begin{center}}
\newcommand{\ece}{\end{center}}
\newcommand{\bpi}{\begin{picture}}
\newcommand{\epi}{\end{picture}}
\newcommand{\bfi}{\begin{figure} \begin{center}}
\newcommand{\efi}{\end{center} \end{figure}}

\newcommand{\bsl}{\begin{slide}{}}
\newcommand{\esl}{\end{slide}}

\newcommand{\pf}{\noindent{\bf Proof}\hspace{7pt}}
\newcommand{\qed}{\rule{1ex}{1ex}}

\newcommand{\hso}[1]{\hspace{-1pt}}

\newcommand{\sbe}{\subseteq}

\newcommand{\sps}{\supset}
\newcommand{\spe}{\supseteq}

\def\<{\langle}
\def\>{\rangle}

\newcommand{\De}{\Delta}

\newcommand{\La}{\Lambda}

\newcommand{\vep}{\varepsilon}

\newcommand{\cA}{{\cal A}}

\newcommand{\cC}{{\cal C}}
\newcommand{\cD}{{\cal D}}

\newcommand{\cF}{{\cal F}}

\newcommand{\cL}{{\cal L}}

\newcommand{\cP}{{\cal P}}
\newcommand{\cQ}{{\cal Q}}

\newcommand{\cS}{{\cal S}}

\newcommand{\Aut}{\mathop{\rm Aut}\nolimits}

\newcommand{\after}{\mathbin{\circ}}

\newcommand{\codim}{\mathop{\rm codim}\nolimits}

\newcommand{\dist}{\mathop{\rm dist}\nolimits}

\newcommand{\id}{\mathop{\rm id}\nolimits}

\newcommand{\mn}{\medskip\noindent}

\def\flexbox#1{\mathchoice{\mbox{#1}}{\mbox{#1}}{\mbox{\scriptsize #1}}%
{\mbox{\tiny #1}}}

\newcommand{\PG}{\mathop{\flexbox{\rm PG}}}

\newcommand{\wA}{{{}^2\hspace{-.2em}A}}
\newcommand{\wD}{{{}^2\hspace{-.2em}D}}

\def\cD{{\cal D}}
\newcommand{\Char}{\mathop{\flexbox{\rm Char}}}

\newcommand{\FF}{{\mathbb F}}

\newcommand{\DP}{{\Delta\hspace{-.1em}\Pi}}

\newcommand{\dfn}{\em}

\newcommand{\Sp}{\mathop{\rm Sp}}

\newcommand{\AI}[1]{\item[\rm{(#1)}]}

\newtheorem{thm}{Theorem}[section]
\newtheorem{prop}[thm]{Proposition}
\newtheorem{cor}[thm]{Corollary}
\newtheorem{lem}[thm]{Lemma}
\newtheorem{conj}[thm]{Conjecture}
\newtheorem{exa}[thm]{Example}
\newtheorem{rem}[thm]{Remark}

\newtheorem{mainthm}{Theorem}

\renewcommand{\qed}{\hfill $\square$}

\newcommand{\Definition}{\refstepcounter{thm}\ \newline\noindent{\bf Definition \thethm}\quad}

\newcommand{\Note}{\refstepcounter{thm}\ \newline\noindent{\bf Note \thethm}\quad}

\newcounter{romanlistctr}
{\end{list}}%

\begin{document}
\title{Projective Subgrassmannians of Polar Grassmannians}
\author{Rieuwert J. Blok \\
Department of Mathematics and Statistics\\
Bowling Green State University\\
Bowling Green, OH 43403\\
U.S.A.\\
{\tt blokr@member.ams.org}\\[5pt]
\\
Bruce N. Cooperstein\\
Department of Mathematics\\
University of California\\
Santa Cruz\\
U.S.A.\\
{\tt coop@ucsc.edu}\\[5pt]}

\date{November 30, 2008
    \begin{flushleft}
        Key Words: Grassmannian, polar geometry, embeddings\\[1em]
    AMS subject classification (2000):
    Primary 51A50;
    Secondary 51E24.
    \end{flushleft}
       }
\maketitle



\newpage

    \begin{abstract}
In this short note, completing a sequence of studies~\cite{CKS,Co2005}, we consider the $k$-Grassmannians of a number of polar geometries of finite rank $n$.
We classify those subspaces that are isomorphic to the $j$-Grassmannian of a projective $m$-space.
In almost all cases, these are parabolic, that is, they are the residues of a flag of the polar geometry.
Exceptions only occur when the subspace is isomorphic to the Grassmannian of $2$-spaces in a projective $m$-space and we describe these in some detail.
This Witt-type result implies that automorphisms of the Grassmannian are almost always induced by automorphisms of the underlying polar space.
\end{abstract}
\newpage

\section{Introduction and preliminaries}\label{sec:introduction}
In \cite{CKS} the authors considered the $k$-Grassmannian $\Gamma$ of a projective $n$-space and asked whether a
 subspace $S$ that is isomorphic to the $l$-Grassmannian of a projective $m$-space is necessarily parabolic.
That is, is it recognizable in the diagram?
They show that this is true in most cases, so that any two isomorphic Grassmannian subspaces
  are conjugate under the action of the automorphism group of $\Gamma$.
In this sense their result is akin to Witt's theorem.
In \cite{Co2005} the second author employs the above result to study a similar question, but now letting $\Gamma$ be a symplectic $k$-Grassmannian.
In the present note, we finish off this sequence of studies by considering the case where $\Gamma$ is the $k$-Grassmannian of an almost arbitrary polar space of finite rank  (See Table~\ref{table:polar types} for the precise list). Our main result, Theorem~\ref{mainthm1}, shows that in contrast to what happens in the cases studied above, not all Grassmannian subspaces are parabolic, although the exceptions are limited.

\mn
For convenience of the reader, before describing our result in more detail, we shall review some basic definitions in Subsection~\ref{subsec:basic}, and we describe the geometries we shall deal with in some detail in Subsection~\ref{subsec:grassmannians}.
The experienced reader should probably skip to Subsection~\ref{subsec:mainthm} for the precise statement of the main result.

\subsection{Basic definitions}\label{subsec:basic}
We assume the reader is familiar with the concepts of a
{\it partial linear rank two
incidence geometry}~\cite{B} $\Gamma = (\Pi, \La)$ (also called a point-line geometry) and the Lie
incidence geometries~\cite{Co} also known as shadow spaces~\cite{Bl,Co1995}.

\mn
The {\it collinearity graph} of $\Gamma$ is the graph $(\Pi, \gamma)$ where
 $\gamma$ consists of all pairs of points belonging to a common line.
For a point $x \in \Pi$ we will denote by $\gamma(x)$ the collection of all points collinear with $x.$
For points $x,y \in \Pi$ and a positive integer $t$ a {\it path of length} $t$
from $x$ to $y$ is a sequence $x_0 = x, x_1, \dots, x_t = y$ such that $\{x_i, x_{i+1}\} \in \gamma$
for each $i = 0, 1, \dots, t-1.$
The {\it distance} from $x$ to $y$, denoted by $d_\Gamma(x,y)$, or simply $d(x,y)$, is defined to be the length of a shortest path from
$x$ to $y$ if some path exists and otherwise is $+\infty.$

\mn By a {\it subspace} of $\Gamma$ we mean a subset $S\subset \Pi$ such that if $l \in \La$ and $l \cap S$ contains at
least  two points, then $l \subset S$.
Clearly the intersection of subspaces is a subspace and consequently it is natural to define the
subspace generated by a subset $X$ of $\Pi, \langle X \rangle_\Gamma,$ to be the intersection of
all subspaces of $\Gamma$ which contain $X$.
\subsection{The geometries under study}\label{subsec:grassmannians}
\paragraph{Projective Grassmannians}
We now describe the $j$-Grassmannian of a projective $m$-space over a field $\FF$.
This is the Lie incidence geometry whose isomorphism type will be denoted $A_{m,j}(\FF)$.

Let  $V$ be a vector space of dimension $m+1$ over $\FF$.
The projective geometry $A_m(\FF)$ is the incidence geometry whose $i$-objects are the (linear) $i$-spaces of $V$, for $i=1,2,\ldots,m$, and in which two objects are incident if one is contained in the other as a subspace.

\mn
For $1 \leq j \leq m,$ let $L_j(V)$ be the collection of all linear $j$-spaces of $V$.
For pairs $(C,A)$ of incident subspaces of $V$ with $\dim(A) = a< j <\dim(C) = c$ let
$$S_j(C,A)=\{B\in L_j(V)\mid A\sbe B\sbe C\}.$$
If $j$ is clear from the context, we shall drop it from the notation.

\mn
The $j$-Grassmannian of $V$, denoted $A_j(V)$ or simply $\cA=(\cP,\cL)$ is the point-line geometry whose point set is $\cP=L_j(V)$ and whose lines are the sets
  $S_j(C,A)$, where $\dim(A) = j-1$ and $\dim(C) = j+1$.
Given a point $p\in \cA$ the collection of points collinear to $p$ is denoted $\alpha(p)$.
\mn

\paragraph{The Polar Grassmannians}
We shall study the Lie incidence geometries that are polar Grassmannians of the types $M_{n,k}(\FF)$ listed in Table~\ref{table:polar types}.
The polar Grassmannian of type $M_{n,k}(\FF)$ is constructed from a non-degenerate reflexive sesquilinear form $\beta$ of Witt index $n$ on a vector space $W$ of dimension $m$ over the field $\FF$ (see ``Construction of a polar Grassmannian'' below).
The type of $\beta$ is given in the table and $m$ is the subscript of the group, which is the full linear isometry group of $\beta$.
In case $\beta$ is $\sigma$-hermitian, we restrict to the case where $\sigma\in\Aut(\FF)$ has order $2$, $\FF$ is a quadratic extension over the fixed field $\FF^\sigma=\{x\in\FF\mid \lambda^\sigma=\lambda\}$, and the norm $N_\sigma\colon\FF\to\FF^\sigma$ is surjective.

\mn
\begin{table}[h]
$$
\begin{array}{@{}c@{}|@{}c@{}|@{}l@{}|c|c|c@{}}
M_{n,k} & \Char(\FF) & \mbox{ $\beta$ }   & \mbox{group} & n & k \\
\hline
B_{n,k}(\FF)& \mbox{ any } & \mbox{ parabolic orthogonal }               & O_{2n+1}(\FF) &\ge 3& 1\le k\le n\\
\hline
C_{n,k}(\FF) & \ne 2 &\mbox{ symplectic }  & \Sp_{2n}(\FF) & \ge 3 & 1\le k\le n \\
\hline
D_{n,k}(\FF) & \mbox{ any } &\mbox{ hyperbolic orthogonal }  & O_{2n}^+(\FF) & \ge 3 & 1\le k\le n-2 \\
\hline
\wA_{2n,k}(\FF) & \mbox{any} &\mbox{ $\sigma$-hermitian }&  U_{2n}^+(\FF) & \ge 3 & 1\le k\le n \\
\hline
\wA_{2n+1,k}(\FF) &\mbox{any} &\mbox{ $\sigma$-hermitian } & U_{2n+1}^+(\FF)& \ge 3 & 1\le k\le n\\
\hline
\wD_{n+1,k}(\FF)  & \mbox{ any }& \mbox{ elliptic orthogonal }& O_{2n+2}^-(\FF)& \ge 3 & 1\le k\le n \\
\end{array}$$
\caption{The polar Grassmannians considered in this paper}\label{table:polar types}
\end{table}

\mn
We shall use the following terminology to distinguish the essentially different geometries for our purposes. Let $M_{n,k}$ be as in Table~\ref{table:polar types}.
If $k=1$, then we call $\Gamma$ a {\dfn polar space}.
In the remaining cases where $k=n$, we call $\Gamma$ a {\dfn dual polar space}.
The remaining polar Grassmannians will be called {\dfn proper polar Grassmannians}.
Occasionally we may call any of these {\dfn orthogonal} if they derive from
$B_n$, $D_n$, or $\wD_{n+1}$ and {\dfn non-orthogonal} otherwise.

\paragraph{Construction of polar Grassmannians}
The polar Grassmannian of isomorphism type $M_{n,k}(\FF)$ is a point-line geometry denoted
 $\Gamma=(\Pi,\La)$. Below we describe the point-set $\Pi$ and line-set $\La$ for each type $M_{n,k}(\FF)$.
Given a point $p\in \Gamma$ the collection of points collinear to $p$ (including $p$) is denoted $\gamma(p)$.

\mn
For a subset $X$ of $W$, let $X^\perp =\{w \in W: \beta(w,x)=0, \forall x \in X\}.$
Let $U$ be a subspace of $W$. Then, $U$ is {\dfn totally isotropic} if $U\sbe U^\perp$.
Recall that the Witt index of $W$ with respect to $\beta$ is $n$.
The non-degenerate {\dfn polar building} of rank $n$ is the incidence geometry $\De$,  whose $i$-objects are the
totally isotropic $i$-spaces of $W$ and in which two objects are incident whenever one contains the other as a subspace.

\mn
For $1\le k\le n$, let $I_k(W)$ be the collection of all $k$-objects of $\De$.
For pairs $(C,A)$ of subspaces of $W$ with $\dim(A) = a< k <\dim(C) = c$ let
$$T_k(C,A)=\{B\in I_k\mid A \subset B \subset C\}.$$
If $k$ is clear from the context, we shall drop it from the notation.
Note that $B\in T_k(C,A)$ forces $A$ to be an object incident to $B$, but $C$ is not necessarily an object of $\De$.
Note also that if $C$ and $A$ are incident objects of $\De$, then $S_k(C,A)=T_k(C,A)$.

\mn
Fix $k$ with $1\le k\le n$ unless $M_n=D_n$ in which case we assume $1\le k\le n-3$.
Then, the {\dfn polar $k$-Grassmannian} $\Gamma=(\Pi,\La)$ of $\De$ is the point-line geometry whose point set is
$\Pi=I_k(W)$ and whose lines are the sets $T(C,A)$, where $A$ is a $(k-1)$-object of $\De$ and $C$ is a $(k+1)$-space of
$W$ with $A\sbe C\sbe A^\perp$.
Note that $k=1$ forces $C$ to be a $2$-object and $A=0$. Similarly, $k=n$ forces $A$ to be an $(n-1)$-object and
$C=A^\perp$.

Finally we construct the geometry of type $D_{n,n-2}$.
We first describe the building $\De$ of type $D_{n}$.
Its objects $X$ are the totally singular (t.s.) subspaces of dimension $1,2,\ldots,n-2,n$.
The type of $X$ is $\dim(X)$ if $1\le i\le n-2$.
The type of $X$ is $n$ or $n-1$ if $\dim(X)=n$.
Another $n$-space $Y$  has the same type as $X$ if $\codim_X(X\cap Y)=\codim_Y(X\cap Y)$ is even.
Incidence is given by inclusion except that two t.s.\ $n$-spaces $Y$ and $Y$ are incident if
 $\codim_X(X\cap Y)=\codim_Y(Y\cap X)=1$.

The point set $\cP$ of $\Gamma$ consists of objects of type $n-2$ and lines are collections of points of the form $$L(L_0,L_-,L_+)=\{P\in\cP\mid L_0\sbe P\sbe L_+\cap L_-\},$$
where $(L_0,L_-,L_+)$ is a flag of type $(n-3,n-1,n)$.

\subsection{Parabolic Grassmannian subspaces of polar Grassmannians.}\label{subsec:mainthm}

\medskip

\noindent When $E \subset F \subset W,$ $\dim(E) = e, \dim(F) = f$ satisfy
$e < k-1, f > k+1$ with $E,F$ totally isotropic,  the collection $T(E,F)$
is a subspace of $\Gamma$ and is isomorphic to an ordinary
Grassmannian geometry $A_{f-e-1,k-e}(\FF).$ Such a subspace is called
\lq\lq parabolic" since the stablizer in $\Aut(\Gamma)$ is a parabolic
subgroup of $\Aut(\Gamma).$  It is natural to ask: Is every subspace of
the polar Grassmannian of type $M_{n,k}(\FF)$ that is isomorphic to some $A_{m,j}(\FF)$ parabolic?

\bmth\label{mainthm1}
Let $\Gamma$ be of type $M_{n,k}(\FF)$ as in Table~\ref{table:polar types} and let
 $S  \cong A_{m,j}(\FF)$ be a subspace of $\Gamma$.
\bit
\AI{i} If $M_n$ is of type $C_n$, $\wA_{2n}$, $\wA_{2n+1}$, or $D_{n,n-2}$, then
 $S$ is parabolic.
\AI{ii} If $M_{n,k}$ is of type $B_{n,k}$, $\wD_{n+1,k}$, or $D_{n,l}$ with $1\le k\le n-3$, then $S$ is parabolic or $k\le n-3$ and $S\cong A_{3,2}(\FF)$ is embedded naturally
 as $D_{3,1}(\FF)$ into the polar subspace $T(C^\perp,C)$, for some $(k-1)$-object $C$.
\eit
\emth
\medskip
\noindent
{\bf Notes on Theorem~\ref{mainthm1}}:
The polar spaces of type $C_n$ , $\wA_{2n-1}$, and $\wA_{2n}$ are similar, but they differ significantly from orthogonal polar spaces in that their hyperbolic lines are longer. This largely explains the division into cases (i) and (ii) of the theorem. The orthogonal $D_{n,n-2}$ geometry appears in case (i) because all exceptions in (ii) are coming from the isomorphism $A_{3,2}\cong D_{3,1}$ which in $D_{n,n-2}$ is merely an isomorphism between two parabolic subspaces.

\mn
The symplectic and unitary cases can be handled in a way similar to~\cite{Co2005}. This has allowed us to leave out a significant portion of the proofs.
Additional arguments are needed to handle the orthogonal Grassmannians. Finally the $D_{n,n-2}$ geometry needs an approach of its own because of the special positioning of node $n-2$ in the diagram.

\medskip
\subsection{Notation and a lemma on convexity of residues}\label{subsec:notation convexity}
Before proceeding to the proofs we introduce some notation:
Since we will generate all kinds of subspaces, of the vector space $W$, of the polar Grassmannian $\Gamma=(\Pi,\La),$ etc. we need to distinguish between these.
When $X$ is some collection of subspaces or vectors from $V$ we will denote the linear subspace of $V$ spanned by $X$ by $\langle X \rangle_V$.
Analogously, we denote the subspace spanned by $X\sbe W$ by $\langle X\rangle_W$.
When $X$ is a set of points in $\Gamma$, we will denote the subspace of $\Gamma$  generated by $X$ by $\langle X \rangle_\Gamma$.
And, when $X$ is a set of points of the projective Grassmannian $\cA$ we will denote the subspace of $\cA$ generated by $X$ by $\langle X \rangle_\cA$.

We will also have to compare distances in $\Gamma$ to distances in some subspace $S$ of $\Gamma$.
The following observation is trivial, but useful.
\ble\label{lem:subspace distances}
Let $\Gamma$ be a partial linear space with subspace $S$ and let $x,y\in S$.
\begin{itemize}
\AI{a} If $d_S(x,y)\le 2$, then $d_\Gamma(x,y)=d_S(x,y)$.
\AI{b} If $d_S(x,y)\ge 3$, then $d_\Gamma(x,y)\le d_S(x,y)$.
Moreover, if $S$ is convex in $\Gamma$, then equality holds.
\end{itemize}
 \ele

\mn
Each polar or projective Grassmannian is the shadow space of some building $\De$. A parabolic subspace is then the shadow of
 a residue $R$ of $\De$. As is probably well-known such a subspace  inherits convexity from the convexity of $R$ in $\De$.
For the benefit of those familiar with the chamber system viewpoint of buildings, we include a proof here.
\ble\label{lem:residues are convex}
If $\Gamma$ is some shadow space of a building $\De$ and $R$ is a residue, then
 the shadow of $R$ in $\Gamma$ is a convex subspace.
\ele
\pf
That the shadow of $R$ in $\Gamma$ is a subspace is trivial. We now show it is convex.
Let $x$ and $y$ be two points on a residue $R$. Let $X$ be a path from $x$ to $y$ containing
 a point $z$ not incident to $R$. Let $\gamma=(c_0,\ldots,c_m)$ be a gallery whose points form $X$ so that
  $c_0\in x\cap R$ and $c_m\in y\cap R$.
Using a retraction $\rho$ onto some apartment containing $c_0$ and $c_m$ we find a gallery supporting a path from $x$ to $y$ whose length is at most that of $X$. Since $\rho$ preserves distances from $c_0$, $\rho(z)$ is still a point outside $R$.
Thus we may assume that $\gamma$ lies in this apartment.
Next we note that  by convexity of $R$ in the chamber system $\De$, $\gamma$ can only be a minimal gallery if it lies within $R$ (and $X$ lies on $R$).
We may then assume that $c_1\not \in R$. Since $\gamma$ is not minimal, the folding ${\mathcal\protect{ f}}$ onto the root $\alpha$ determined by $\{c_0,c_1\}$ and containing $c_0$ sends $\gamma$ to a gallery from $c_0$ to $c_m$. Moreover, the points on $c_0$ and ${\mathcal\protect{f}}(c_1)$ coincide. Therefore the corresponding path in $\Gamma$ is strictly shorter than $X$. It follows that $X$ is not a geodesic in $\Gamma$ and that any geodesic must consist of points in the shadow of $R$.
\qed

\section{Properties of Projective Grassmannians}\label{section:projective grassmannians}
\mn In this short section we recall some properties of a projective Grassmannian incidence geometry $\cA=A_j(V)$ of type $A_{m,j}(\FF)$. We omit the proofs because they are either well known or entirely straightforward to prove.

\ble\label{lem:up or down}
Suppose $2\le j\le m-1$.
If $x,y,z$ are points of $A_{m,j}$ on distinct lines
 $S(E_i,D_i)$, $i=1,2,3$, then either $E_1=E_2=E_3$ or $D_1=D_2=D_3$.
\ele

\begin{lem}\label{lem:sing in Grass}
\bit
\AI{i}
There are two classes of maximal singular subspaces of $\cA=(\cP,\cL)$ with representatives
$S(V,D)$ where $\dim(D) = j-1$ and $S(E,0)$ where $\dim(E) = j+1$.
Then, $S(V,D) \cong A_{m+1-j,1}(\FF)$ and
$S(E,0) \cong A_{j,j}(\FF).$    Those of the first class will be referred to as type $+$ and the second class as type $-$.

\AI{ii}  If $M_1$ and $M_2$ are maximal singular subspaces and $M_1 \cap M_2$ is a line then $M_1$ and $M_2$ are in different classes.  If $M_1 \cap M_2$ is a point then they are in the same class.
\eit
\end{lem}

\begin{lem}\label{generation by type + maximal singulars}

\noindent Let $M$ be a maximal singular subspace of $\cA$ of type $+$.  Then
$\langle M \rangle_V = V$.

\end{lem}

\medskip

\noindent Now let $U$ be a hyperplane of $V$ and $X$ a  $1$-space of $V$ with $X$ not contained in $U$.  Set
$\cP(U) = \{x \in \cP:x \subset U\}$ and $\cP_X = \{x \in \cP: X \subset x\}$.

\medskip

\begin{lem}\label{decomp of Grass}
\bit
\AI{i} $\cP(U)$ is a subspace of $\cP$ and $\cP(U) \cong A_{m-1,j}(\FF).$

\AI{ii} $\cP_X$ is a subspace of $\cP$ and $\cP_X \cong A_{m-1,j-1}(\FF).$

\AI{iii} If $x \in \cP(U)$ then $\alpha(x) \cap \cP_X$ is a maximal singular subspace of type $-$ in $\cP_X$ isomorphic to $A_{j-1,j-1}(\FF)$.
Furthermore, $\langle x, \alpha(x) \cap \cP_X \rangle_\cA$ is a maximal singular subspace of type $-$ in $\cP$.

\AI{iv} If $y \in \cP_X$ then $\alpha(y) \cap \cP(U)$ is a maximal singular subspace of type $+$ in $\cP(U)$ isomorphic to $A_{m-j,1}(\FF)$.
Futhermore, $\langle y, \alpha(y) \cap \cP(U)\rangle_\cA$ is a maximal singular subspace of type $+$ in $\cP$.

\AI{v} If $x_1,x_2 \in \cP(U)$  are collinear then $\alpha(x_1) \cap \alpha(x_2) \cap \cP_X$ is a point.
Similarly, if $y_1,y_2 \in \cP_X$ are collinear then $\alpha(y_1) \cap \alpha(y_2) \cap \cP(U)$ is a point.
\eit
\end{lem}

\begin{lem}\label{diameter of Grassmannian}

\noindent The diameter of the collinearity graph of $A_j(V)$ is $min\{j,m+1-j\}$. For $x,y \in \cP, d(x,y) = \dim(x/x \cap y) = \dim(y/x \cap y)$.

\end{lem}

\section{Properties of polar Grassmannians}\label{section:polar grassmannians}
In this section we review some properties of polar Grassmannians of type $M_{n,k}(\FF)$ as listed in Table~\ref{table:polar types}.

We first study the singular subspaces of $\Gamma$.
In the interest of brevity, just like in Section~\ref{section:projective grassmannians}, we omit the proofs because these are either well known or easy to reproduce.

\begin{lem}\label{lem:sing in sing}
Let $\Gamma$ be of type $M_{n,k}$ as in Table~\ref{table:polar types}. Then $\Gamma$ is a subgeometry of $A_k(W)$. Hence, the singular subspaces of $\Gamma$ are contained in
 singular subspaces of $A_k(W)$.\qed
\end{lem}

\mn 
Lemmas~\ref{lem:sing in sing}~and~\ref{lem:sing in Grass} imply the following lemma.
\newpage
\begin{lem} \label{lem:sing in proper}\
\vspace{.1em}
Suppose $M_{n,k}(\FF)$ is as in Table~\ref{table:polar types}.
\bit
\AI{i}The polar Grassmannian space $\Gamma$ has two classes of maximal singular subspaces with representatives $T(B,0)$ where $B$ is a totally isotropic subspace of $W$, with $\dim(B) = k + 1,$ and $T(C,A)$ where $A$ and  $C$ are incident  totally isotropic subspaces of $W$, where  $\dim(A) = k-1$, and $\dim(C) = n$.
In the former case $T(B,0) \cong A_{k,k}(\FF)$ and in the latter $T(C,A) \cong
A_{n-k,1}(\FF)$.  We refer to the first as type $-$ maximal singular subspaces and the latter as type $+$.

\AI{ii}
If $M_1$ and $M_2$ are maximal singular subspaces  of $\Gamma$ of different types, then $M_1 \cap M_2$ is either empty or a line.

\AI{iii} 
If $M_1$ and $M_2$ are distict type $-$ maximal singular subspaces of $\Gamma$, then $M_1 \cap M_2$ is either empty or a point.
\eit
\end{lem}

\Note If in part (i) of Lemma~\ref{lem:sing in proper} we have  $M_{n,k}=D_{n,k}$, then $C$ can be of type $n-1$ or $n$.
In addition, in cases $M_{n,1}$, $D_{4,3}$, $D_{4,4}$ -- the latter two viewed as having type $D_{4,1}$ -- the maximal singular subspace of type $+$ has $A=\{0\}$.

\Definition
\noindent A {\dfn symp} of $\Gamma$ or $\cA$ is a maximal geodesically closed subspace which is isomorphic to a non-degenerate polar space.

\mn We shall now determine all symps of $\Gamma$.
Note that every non-degenerate polar space is the convex closure of any two of its points at distance $2$ from each other.
Therefore we determine such pairs of points in $\Gamma$ and describe their convex closure. If this convex closure is {\em not}
 a symp, we call the pair of points {\dfn special}.

\ble\label{lem:points at distance two}
\bit\AI{i}
Let $M_{n,k}$ be one of $B_{n,1}$, $C_{n,1}$, $\wA_{2n,1}$, $\wA_{2n+1,1}$, $\wD_{n+1,1}$, $D_{n,1}$ with $n\ge 4$.
Then there is one class of points at distance two in $\Gamma$.  For one such pair $(x,y)$ as subspaces of $W$, $\dim(x \cap y) = 0$. The unique symp on $\{x,y\}$ is $\Gamma$ itself.

\AI{ii}
Let $M_{n,k}$ be one of $B_{n,k}$, $C_{n,k}$, $\wA_{2n,k}$, $\wA_{2n+1,k}$, $\wD_{n+1,k}$, with $2\le k\le n-1$ or be $D_{n,k}$ with $2\le k\le n-3$.
Then, there are three classes, of type $+$, $0$, and $-$, of points at distance two in $\Gamma$.

A pair $(x,y)$ is of $-$ type if, as subspaces of $W$, $\dim(x \cap y) = k-2$ and $x \perp y$.
The unique symp of $-$ type on $\{x,y\}$ is $T(x + y, x \cap y)\cong A_{3,2}$.

A pair $(x,y)$ is of $+$ type if, as subspaces of $W$, $\dim(x \cap y) = k-1$ and $(x + y)/(x \cap y)$ is a non-degenerate $2$-space.
The unique symp of $+$ type on such a pair is $T((x \cap y)^\perp, x \cap y)\cong M_{n+1-k,1}$.

A pair $(x,y)$ is of $0$ type if, as subspaces of $W$, $\dim(x\cap y)=k-2$, and $z=(x + y)\cap (x+y)^\perp$ is a  point of $\Gamma$.
Now $(x,y)$ is called a {\em special pair}, the convex closure of $\{x,y\}$ is the set of points on the two lines $xy$ and $yz$, which is not a symp.

\AI{iii}
Let $M_{n,k}$ be one of $B_{n,n}$, $C_{n,n}$, $\wA_{2n,n}$, $\wA_{2n+1,n}$, $\wD_{n+1,n}$.
Then there is one class of points at distance two in $\Gamma$.  For one such pair $(x,y)$, as subspaces of
$W$, $\dim(x \cap y) = n-2$. The unique symp on such a pair is $T((x \cap y)^\perp, x \cap y)\cong M_{2,2}$.

\eit
\ele

\bco\label{cor:points in symps}
Let $M_{n,k}$ be as in Table~\ref{table:polar types}.
If $x$ and $y$ are two points at distance two from each other in $\Gamma$ having at least two common neighbors, then
 they are contained in a unique symp as described in Lemma~\ref{lem:points at distance two}.
\eco
\pf
The condition on the number of common neighbors excludes the possibility that $x$ and $y$ form a special pair.
The result follows from Lemmma~\ref{lem:points at distance two}.
\qed

\section{Proof of the Main Theorem}\label{section:almost all cases}

\medskip

\noindent In this section we handle all cases of Theorem~\ref{mainthm1} with the exception of the $D_{n,n-2}$-geometry, which requires special consideration and will be dealt with in Section~\ref{section:Dn,n-2}.
\medskip

\subsection{Polar spaces and dual polar spaces}
In this first lemma, we encounter the exceptional $A_{3,2}$-subspaces
 of Theorem~\ref{mainthm1} (ii).
\ble\label{lem:polar j=2}
\bit
\AI{i}
Let $\Gamma$ be a polar space of type $C_{n,1}$ with $\Char(\FF)\ne 2$, or one of type $\wA_{2n,1}$, $\wA_{2n+1,1}$.
If $\Gamma$ contains a subspace $S\cong A_{m,j}$, then $1\le m\le n$, $j=1$ and $S$ is parabolic.

\AI{ii}
Let $\Gamma$ be a polar space of type $B_{n,1}$, $D_{n,1}$ with $n\ge 3$, or  $\wD_{n+1,1}$.
If $\Gamma$ contains a subspace $S\cong A_{m,j}$, where $\min\{j,m+1-j\} = 2$ and $m\ge 3$, then in fact
  $S\cong A_{3,2}$ and it is embedded as $D_{3,1}$ into a non-degenerate $6$-space of the natural embedding $W$ of
   $\Gamma$.
\eit
\ele
\pf
Obviously subspaces of type $A_{m,1}$ exist and must be parabolic.
If $\Gamma$ contains a subspace of type $A_{m,j}$ with $\min\{j,m+1-j\}\ge 2$, then it also contains a subspace $S'\cong A_{3,2}$.
Let $U = \langle S'\rangle_W$ a vector subspace of $W$.
Since $S'$ is strongly hyperbolic (see [CS]) it follows that $\dim(U) = 6$.
We now have an embedding from the orthogonal polar space $S'$ into the polar space $G\cong M_{n,1}$.
Consider a polar frame $\cF=\{e_1,f_1,e_2,f_2,e_3,f_3\}$ in $S'$.
By Blok and Brouwer~\cite{BB} and Cooperstein and Shult~\cite{CoSh} we know that this set generates $S'$.

(i)
In the cases $\wA_{2n,k}$ and $\wA_{2n+1,k}$, $\cF$ generates a subspace of type $\wA_{6,1}\not\cong S'$, a contradiction (See~\cite{BlCo}).
In the case $C_{n,k}$, since $\Char(\FF)\ne 2$, $\cF$ generates a subspace of type $C_{3,1}\not\cong S'$, again a contradiction. This concludes case (i).

(ii)
In the cases $B_{n,1}$, $D_{n,1}$ with $n\ge 3$, or  $\wD_{n+1,1}$, $\cF$ generates the full subspace $T(U,0)$ of type $D_{3,1}$ supported by $U$.
Clearly $S'\sbe T(U,0)$.
Since subspaces of type $D_{3,1}(\FF)$ do not contain proper subspaces isomorphic to themselves we have
 $S'=T(U,0)$.

We now show that in case (ii) $S$ cannot be of type $A_{m,j}$ with $\min\{j,m+1-j\} = 2$ and $m\ge 4$.
Suppose there is such a subspace. Then there is also a subspace $S'\cong A_{4,2}$.
Then $S'$ contains a point $p$ and a line $L$ no point of which is collinear to $p$.
This contradicts the fact that the Buekenhout-Shult axioms hold in the polar space $\Gamma$.
\qed

\mn
It is easy to see that dual polar spaces don't have proper parabolic projective Grassmannian subspaces.
For completeness we include the following lemma.
\begin{lem}\label{lem:A32 dual polar}
The dual polar spaces of type $B_{n,n}$, $C_{n,n}$, $\wA_{2n,n}$, $\wA_{2n+1,n}$, and $\wD_{n+1,n}$ have no subspaces
 of type $A_{m,j}(\FF)$ such that $\min\{m+1-j,j\}\ge 2$.
\end{lem}

\medskip

\noindent {\bf Proof}:
If on the contrary there is such a subspace, then there is also a subspace $S \cong A_{3,2}(\FF)$.
Since $S$ is a polar space it is contained in some symp $\cS$ of $\Gamma$ by Corollary~\ref{cor:points in symps}.
By Lemma~\ref{lem:points at distance two}, the only type of symp is of the form $T(C^\perp,C)$, where $\dim(C)=n-2$.
Now $\cS$ has polar rank $2$, whereas $S$ has polar rank $3$, a contradiction.
\qed

\subsection{Proper polar Grassmannians}
We now consider the cases where $\Gamma$ is a proper polar Grassmannian.
Our proof is by induction on $N = n + k + m + \min\{j,m+1-j\}$.

\begin{lem}\label{lem:A32 non-orthogonal subspaces}
Let $\Gamma$ be of type $\wA_{2n,k}$ or $\wA_{2n+1,k}$, or $C_{n,k}(\FF)$ with $\Char(\FF)\ne 2$, and
 let $2\le k\le n-1$.
Then any subspace $S \cong A_{3,2}(\FF)$ of $\Gamma$ is parabolic.
\end{lem}

\medskip

\noindent {\bf Proof}:
Assume  $S \cong A_{3,2}(\FF)$. Since $S$ is a polar space it is contained in some symp $\cS$ of $\Gamma$
By Corollary~\ref{cor:points in symps}.
By Lemma~\ref{lem:points at distance two} there are two possibilities for $\cS$:

($-$) If $\cS$ is of $-$ type, then there are totally isotropic subspaces $D \subset E$ such that $\dim(D) = k-2$, $\dim(E) = k + 2$ with $\cS = T(E,D)$. In this case we have $S \cong T(E,D)$.
Since a geometry of type $A_{3,2}$ does not have proper subspaces isomorphic to itself, we find $S=\cS$, which is parabolic.

($+$) If $\cS$ is of $+$ type, then there is a totally isotropic subspace $C$, $\dim(C) = k-1$ such that $\cS = T(C^\perp,C)$.
Let $U = \langle S\rangle_W$ a vector subspace of $C^\perp$.
The map taking $x \in S$ to $x/C$ is an embedding of the polar space $S$ into $\PG(U/C)$.
Since $S$ is strongly hyperbolic (see [CS]) it follows that $\dim(U/C) = 6$.
We now have an embedding from the orthogonal polar space $S$ into the polar space $T(C^\perp,C)$ of type $M_{n+1-k,1}$.
In the cases $\wA_{2n,k}$ and $\wA_{2n+1,k}$, this is impossible because the smallest subspace of $\wA_{2n,1}$
 spanned by three hyperbolic lines is $\wA_{6,1}$, which is not isomorphic to $D_{3,1}$.
In the case $C_{n,k}$, since $\Char(\FF)\ne 2$, the smallest subspace of $C_{n,1}$ containing three pairwise orthogonal hyperbolic lines is $C_{3,1}$, which is not isomorphic to $D_{3,1}$. This contradiction concludes the proof.
\qed

\begin{lem}\label{lem:A32 orthogonal subspaces}
Let $\Gamma$ be of type $B_{n,k}$ or $\wD_{n+1,k}$ with $2\le k\le n-1$,  or $D_{n,k}$ with $2\le k\le  n-3$.
Then a subspace $S \cong A_{3,2}(\FF)$ of $\Gamma$ is parabolic, or it is embedded
as $T(U,C)\cong D_{3,1}$, where $\dim(C)=k-1$, $\dim(U)=k+5$, and $U/C$ is non-degenerate.
\end{lem}

\medskip

\noindent {\bf Proof}:
The proof is almost identical to that of Lemma~\ref{lem:A32 non-orthogonal subspaces} with the following exception.
In case $\cS$ is a symp of $+$ type, as
 in Lemma~\ref{lem:polar j=2}, $S$ is embedded into $T(C^\perp,C)$ of type $M_{n+1-k,1}$ as $T(U,C)$, where $\dim(C)=k-1$,
 $\dim(U)=k+5$, and $U/C$ is non-degenerate.
\qed

\ble\label{lem:regular j=2 m ge 4}
Let $\Gamma$ be of type $B_{n,k}$, $\wD_{n+1,k}$, $C_{n,k}$, $\wA_{2n,k}$, $\wA_{2n+1,k}$ and $2\le k\le n-1$,
 or $D_{n,k}$ with $2\le k\le n-3$.

If $S\cong A_{m,j}$ is a subspace of $\Gamma$, where $\min\{j,m+1-j\} = 2$ and $m\ge 4$, then $S$ is parabolic.
\ele

\pf
By assumption $S$ has a subspace $S'\cong A_{4,2}(\FF)$.
Let $\cD$ be a symp of $S'$.
Since $\cD\cong A_{3,2}$ is a polar space it is contained in a symp $\cS$ of $\Gamma$ by Corollary~\ref{cor:points in symps}.
According to Lemma~\ref{lem:points at distance two} there are two types of symps.

Either $\cS$ is of $+$ type and we have $\cS=T(C^\perp,C)$ for some t.i.\ $(k-1)$-space $C$,
 or $\cS$ is of $-$ type and we have $\cS=S(B,A)$, where $A$ is a t.i.\ $(k-2)$-space and
  $B\sps A$ is a t.i.\ $(k+2)$-space.

We claim that $\cS$ cannot be of $+$ type.
Namely, we claim that $S'\sbe \cS$ contradicting Lemma~\ref{lem:polar j=2}.
To see this, suppose that $x\in S'-\cD$.
Then, considering the $A_{4,2}$ geometry $S'$ we see that $\gamma(x)\cap \cD$ is a projective plane.
In particular, $x$ is collinear to three pairwise collinear points, not all on one line and all containing
 $C$.
By Lemma~\ref{lem:up or down}, $C\sbe x$ so that $x\in \cS$.
This proves the claim.

\mn
This contradiction shows that $\cS=S(B,A)$ where $A$ is a t.i.\ $(k-2)$-space and
  $B\sps A$ is a t.i.\ $(k+2)$-space.
{}From this it follows that if $x,y \in S$ with $d(x,y) = 2$, then as subspaces of $W$ they satisfy $x \perp y$
by Lemma~\ref{lem:points at distance two}.
Since the diameter of $S$ is two by Lemma~\ref{diameter of Grassmannian} it then follows that $B = \langle S \rangle_W$ is a totally isotropic subspace of $W$.
Consequently, $S \subset T(B,0)$.  By Theorem (2.15) of \cite{CKS} it follows that $S$ is parabolic and the theorem holds.
\qed

\mn {\bf Remark:} We'd like to point out that the proof of Lemma~\ref{lem:regular j=2 m ge 4} does not rely on
Lemmas~\ref{lem:A32 non-orthogonal subspaces} and \ref{lem:A32 orthogonal subspaces}.

\mn Having proved Lemmas~\ref{lem:A32 non-orthogonal subspaces} and ~\ref{lem:A32 orthogonal subspaces}, we may assume that $m \geq 5$ and $\min\{j,m+1-j\} \geq 3.$
We continue with the notation of the introduction where $V$ was introduced as an $(m+1)$-dimensional vector space and $\cA=(\cP,\cL)$ is the Grassmannian geometry of $j$-dimensional subspaces of $V$.  Let $\tau: \cP \to S$ be an isomorphism of geometries. As in Section~\ref{section:projective grassmannians} let $U$ be a hyperplane of $V$ and $X$ a one-dimensional subspace of $V$ such that $X$ is not contained in $U$
and set $\cP(U) = \{x \in \cP: x \subset U\}$ and $\cP_X = \{x \in \cP: X \subset x\}$.
Also, set $ S_U = \tau(\cP(U))$ and $S_X = \tau(\cP_X)$.

\mn
Since $S_U \cong A_{m-1,j}(\FF)$ and $(m-1) + \min\{j, m - j\} < m + \min\{j,m+1-j\}$ it follows by our induction hypothesis that $S_U = T(B_U,A_U)$
where $A_U \subset B_U$ are totally isotropic subspaces with $\dim(A_U) = a_U$, $\dim(B_U) = b_U$ and $m-1 = b_U - a_U$,
$j = k - a_U$.

\mn
Similarly,  since $S_X \cong A_{m-1,j-1}(\FF)$ and $(m-1) + \min\{j-1, m - (j-1)\} < m + \min\{j,m+1-j\}$ it follows that $S_X = T(B_X,A_X)$ where $A_X \subset B_X$ are totally isotropic subspaces with $\dim(A_X) = a_X$,
$\dim(B_X) = b_X$ and $m = b_X - a_X$, $j-1 = k - a_X$.

\mn
Let $x \in S_U, y \in S_X$ be collinear.  Then by Lemma~\ref{decomp of Grass},
$Y_-= \langle x, S_X \cap \gamma(x)\rangle_\Gamma$ and $Y_+ = \langle y, S_U \cap \gamma(y)\rangle_\Gamma$ are maximal singular subspaces of $S$ which meet in a line.
Namely, $\tau^{-1}(Y_-)=S(\tau(x)^{-1}+\tau(y)^{-1},0)$, whereas
 $\tau^{-1}(Y_+)=S(V,\tau(x)^{-1}\cap\tau(y)^{-1})$, which meet in the line
  $\tau^{-1}(Y_-\cap Y_+)=S(\tau(x)^{-1}+\tau(y)^{-1},\tau(x)^{-1}\cap\tau(y)^{-1})$.

\mn
Let $M_\vep$ be a maximal singular subspace of $\Gamma$ containing $Y_\vep$, $\vep\in\{+,-\}$.
By Lemma~\ref{lem:sing in proper}, $M_\vep$ can be of type $+$ or $-$.
Since $M_-$ and $M_+$ are distinct and intersect in at least a line, they are of different type, again by Lemma~\ref{lem:sing in proper}.
Consequently, at least one of $M_-$, $M_+$ is of type $+$.
For the sake of argument, assume $M_+$ is of type $+$.
Then there is a maximal totally isotropic subspace $B$ and a $(k-1)$-dimensional subspace $A \subset B$ such that $M_+ = T(B,A)$.

\mn
Now consider $M_+ \cap S_U =  T(B,A) \cap T(B_U,A_U)$.
Since $Y_+\sbe M_+$ we have $M_+\cap S_U\spe S_U\cap \gamma(y)$.
On the other hand, since $y\in M_+$ and $M_+$ is singular we have $M_+\cap S_U\sbe S_U\cap \gamma(y)$.
Thus, $M_+\cap S_U=S_U\cap \gamma(y)$.
By looking at $\cP(U)$ we see that $S_U\cap\gamma(y)$ is in fact a maximal singular subspace of $S_U=T(B_U,A_U)$.
It follows that $A_U\sbe A$ and $B_U\sbe B$ so that $M_+\cap S_U=T(B_U,A)$.

Then $B_U = \langle M_+ \cap S_U\rangle_W = \langle S_U \cap \gamma(y) \rangle_W$ by
Lemma~\ref{generation by type + maximal singulars}
which implies that $B_U \subset y^\perp$ since $x^\prime \in \gamma(y)$ implies $y \perp x^\prime$.

\mn
Now assume that $y^\prime \in S_X$ such that $y^\prime, y$ are collinear.
Then by looking at $V$, we find that $S_U \cap \gamma(y)$ and $S_U \cap \gamma(y^\prime)$ are
 maximal singular subspaces of $S_U$ intersecting in a point.
 Hence, they must be in the same class of maximal singular subspaces of $S_U$.
Applying the argument for $y$ to $y'$, we find that $B_U \subset (y^\prime)^\perp$.

\mn Since the collinearity graph of $S_X$ is connected, it follows that for all $z \in S_X$, we have $B_U \subset z^\perp$.
Since $\langle S_X \rangle_W = B_X$ we have $B_X \perp B_U$.

\mn
Set $D = B_U + B_X$, a totally isotropic subspace.  Now $S_U, S_X \subset T(D,0).$
Since $\langle S_U,S_X\rangle_\cA = S$ (as follows from~\cite{BB,CoSh,RS}),
it follows that $S \subset T(D,0),$ for, if $x,y$ are collinear points of $\Pi$ and $x,y \subset D$,
then for every $z \in T(x+y, x \cap y)$ also $z \subset D$.
Now  we are done by Theorem (2.15) of~\cite{CKS}.
\qed

\section{The geometry of type $D_{n,n-2}$}\label{section:Dn,n-2}
In this section we consider geometries $\Gamma$ of type $D_{n,n-2}(\FF)$, for some arbitrary field $\FF$.
Our techniques here are significantly different from those in Section~\ref{section:almost all cases}.
In Lemmas~\ref{lem:A3,2 in Dn,n-2}~and~\ref{lem:no A5,3} we consider pairs of points at distance $2$ or $3$.
Lemma~\ref{lem:Dk+2} allows us to determine the convex subspaces containing such pairs.

\ble\label{lem:Dk+2}
If $x$ and $y$ are points of $\Gamma$ at distance $k$ in the collinearity graph, then
 $x$ and $y$ share an object of type $n-2-k$, so $x$ and $y$ are contained in a subgeometry of type $D_{k+2,k}(\FF)$.
Moreover, this subgeometry is a convex subspace of $\Gamma$.
\ele
\pf
The first part follows immediately from the definition of collinearity.
The second part follows from Lemma~\ref{lem:residues are convex}.
\qed

\ble\label{lem:A3,2 in Dn,n-2}
Subspaces of $\Gamma$ isomorphic to $A_{3,2}(\FF)$ are parabolic.
\ele
\pf
Let $S$ be a subspace isomorphic to $A_{3,2}(\FF)$ and let $x,y\in S$ be at distance $2$. Note that  $S$ has diameter $2$ and is the convex closure, in $S$, of any two points at distance $2$ in $S$. By Lemma~\ref{lem:subspace distances}, also $d_\Gamma(x,y)=2$ and $S$ is contained in the convex closure of $\{x,y\}$ in $\Gamma$.
Hence, by Lemma~\ref{lem:Dk+2}, we may assume that $n=4$ and $k=2$.

Considering an apartment of type $D_4$ on $x$ and $y$ one finds the following two possibilities.
(1)  $x$ and $y$ form a "special pair", that is, they have exactly one common neighbor. Since $x$ and $y$ already have many common neighbors in $S$, this is impossible.
(2) $x$ and $y$ are at symplectic distance $2$, that is they share an object $R$ of type $1$, $3$, or $4$.
Note that $R$ and $S$ are isomorphic. Both are the convex closures of the two points in $S$ so that
 $S\sbe R$. Moreover, both $R$ and $S$ are generated by the six points of an apartment, so that we have $R=S$. Thus, $S$ is parabolic.
\qed

\mn
In the proof of Lemma~\ref{lem:no A5,3} we shall study a pair of points $x$ and $y$ at distance $3$ from each other and a collection $\cC$ of points on some geodesic from $x$ to $y$. Lemma~\ref{lem:Dk+2} says that we can study that situation entirely within a geometry of type $D_{5,3}$.
Since all possible configurations of pairs of points are realized inside an apartment of the $D_5$-building, we shall describe this apartment here.
Recall that $\Gamma$ is constructed from the polar space on the $10$-dimensional vector space $W$ over $\FF$ endowed with the non-degenerate symmetric bilinear form $\beta$ such that
 $\cQ(x)=\beta(x,x)$ is a non-degenerate hyperbolic quadratic form on $W$.
Let $\{e_i,f_i\}_{i=1}^5$ be a hyperbolic basis for $W$ with respect to $\beta$, so that $\beta(e_i,f_i)=1=\beta(f_i,e_i)$ for $i=1,2,3,4,5$, and all other inner products between these basis vectors are zero.
The apartment $\Sigma$ for the $D_{5,3}$-geometry
 supported by this basis consists of all points and lines of $\Gamma$ spanned, as subspaces of $W$, by subsets of this basis.

Now fix the point $X=\langle e_1,e_2,e_3\rangle_W$.
In Figure~\ref{fig:distances in D53}, for each distance class relative to $X$, we have listed a representative point in $\Sigma$.
We call the classes $C_i$, where $i\in\{0,1,2g,2q,2s,3h,3q,3hh, 4\}$.
If $x$ and $y$ are points such that $y\in C_i(x)$, then we shall write $\dist_\Gamma(x,y)=i$.
An arrow indicates that collinear representatives exist.

\begin{figure}[h]
$$\begin{array}{lrclrclr}
0 : & \langle e_1,e_2,e_3\rangle_W  & \rightarrow & 1 : & \langle e_1,e_2,e_4\rangle_W  & \rightarrow  &2 g: & \langle e_1,e_4,e_5\rangle_W  \\
   && \swarrow  & \downarrow && \swarrow & \downarrow & \\

2q: &\langle e_1,e_2,f_3\rangle_W & \rightarrow & 2s: &\langle e_1,f_3,e_4\rangle_W & \rightarrow & 3h: &\langle f_3,e_4,e_5\rangle_W\\
   && \swarrow  & \downarrow && \swarrow &  & \\
3 q: & \langle e_1,f_2,f_3\rangle_W & \rightarrow  & 3hh: &\langle f_2,f_3,e_4\rangle_W& \rightarrow  &4 : & \langle f_1,f_2,f_3\rangle_W\\
\end{array}$$
\caption{Representative objects for all distance classes in a $D_{5,3}$ apartment.}\label{fig:distances in D53}
\end{figure}

Analyzing the apartment $\Sigma$, we find the following.
\ble\label{lem:distance distribution of D5,3}
Figure~\ref{fig:distance distribution of D53} is the distance distribution diagram for the $D_{5,3}$-apartment.
\ele
\begin{figure}[h]
\setlength{\unitlength}{1pt}
\bpi(220,210)(-100,0)
\thicklines{
\multiput(10,10)(100,0){3}{\circle{20}}
\multiput(20,10)(100,0){2}{\line(1,0){80}}
\multiput(110,20)(100,0){1}{\line(0,1){80}}
\multiput(16,16)(100,0){2}{\line(1,1){87}}

\multiput(10,110)(100,0){3}{\circle{20}}
\multiput(20,110)(100,0){2}{\line(1,0){80}}
\multiput(110,120)(100,0){2}{\line(0,1){80}}
\multiput(16,116)(100,0){2}{\line(1,1){87}}

\multiput(10,210)(100,0){3}{\circle{20}}
\multiput(20,210)(100,0){2}{\line(1,0){80}}

\put(28,2){\makebox(0,0){$\scriptstyle C_{3q}$}}
\put(128,2){\makebox(0,0){$\scriptstyle C_{3hh}$}}
\put(228,2){\makebox(0,0){$\scriptstyle C_{4}$}}

\put(28,102){\makebox(0,0){$\scriptstyle C_{2q}$}}
\put(128,102){\makebox(0,0){$\scriptstyle C_{2s}$}}
\put(228,102){\makebox(0,0){$\scriptstyle C_{3h}$}}

\put(28,202){\makebox(0,0){$\scriptstyle C_{0}$}}
\put(128,202){\makebox(0,0){$\scriptstyle C_{1}$}}
\put(228,202){\makebox(0,0){$\scriptstyle C_{2g}$}}

\put(10,10){\makebox(0,0){${\scriptstyle 3}$}}
\put(110,10){\makebox(0,0){${\scriptstyle 12}$}}
\put(210,10){\makebox(0,0){${\scriptstyle 1}$}}

\put(10,110){\makebox(0,0){${\scriptstyle 3}$}}
\put(110,110){\makebox(0,0){${\scriptstyle 24}$}}
\put(210,110){\makebox(0,0){${\scriptstyle 12}$}}

\put(10,210){\makebox(0,0){${\scriptstyle 1}$}}
\put(110,210){\makebox(0,0){${\scriptstyle 12}$}}
\put(210,210){\makebox(0,0){${\scriptstyle 12}$}}

\multiput(-5,10)(0,100){3}{\makebox(0,0){${\scriptstyle 0}$}}

\multiput(225,110)(0,100){2}{\makebox(0,0){${\scriptstyle 2}$}}
\put(225,10){\makebox(0,0){${\scriptstyle 0}$}}

\multiput(95,15)(0,100){3}{\makebox(0,0){${\scriptstyle 1}$}}
\multiput(195,115)(0,100){2}{\makebox(0,0){${\scriptstyle 4}$}}

\put(25,215){\makebox(0,0){${\scriptstyle 12}$}}
\put(25,115){\makebox(0,0){${\scriptstyle 8}$}}
\put(25,15){\makebox(0,0){${\scriptstyle 4}$}}

\put(125,215){\makebox(0,0){${\scriptstyle 4}$}}
\put(125,115){\makebox(0,0){${\scriptstyle 2}$}}
\put(125,15){\makebox(0,0){${\scriptstyle 1}$}}

\put(105,195){\makebox(0,0){${\scriptstyle 2}$}}
\put(105,125){\makebox(0,0){${\scriptstyle 1}$}}
\put(105,95){\makebox(0,0){${\scriptstyle 1}$}}
\put(105,25){\makebox(0,0){${\scriptstyle 2}$}}

\put(205,195){\makebox(0,0){${\scriptstyle 2}$}}
\put(205,125){\makebox(0,0){${\scriptstyle 2}$}}

\put(125,95){\makebox(0,0){${\scriptstyle 2}$}}

\put(195,15){\makebox(0,0){${\scriptstyle 12}$}}

\put(17,124){\makebox(0,0){${\scriptstyle 4}$}}
\put(17,24){\makebox(0,0){${\scriptstyle 8}$}}
\put(117,124){\makebox(0,0){${\scriptstyle 2}$}}
\put(117,24){\makebox(0,0){${\scriptstyle 4}$}}

\put(96,203){\makebox(0,0){${\scriptstyle 1}$}}
\put(96,103){\makebox(0,0){${\scriptstyle 1}$}}
\put(196,203){\makebox(0,0){${\scriptstyle 4}$}}
\put(196,103){\makebox(0,0){${\scriptstyle 4}$}}

\put(110,-5){\makebox(0,0){${\scriptstyle 4}$}}
\put(110,225){\makebox(0,0){${\scriptstyle 4}$}}
}\epi
\caption{The distance distribution diagram of the $D_{5,3}$ apartment.}\label{fig:distance distribution of D53}
\end{figure}
We briefly explain the meaning of Figure~\ref{fig:distance distribution of D53}.
The classes $C_i$, where $i\in\{0$, $1$, $2g$, $2q$, $2s$, $3h$, $3q$, $3hh$, $4\}$, are represented by balloons and are arranged exactly as in Figure~\ref{fig:distances in D53}.
The numbers inside the balloons indicate the size of the class.
The edge between classes $C_i$ and $C_j$ has two labels: $n_{ij}$ and $n_{ji}$. The label $n_{ij}$ close to $C_i$ indicates the number of neighbors in $C_j$ for each member of $C_i$.
Each class $C_i$ also has a label $n_{ii}$ indicating the number of neighbors in $C_i$ for each member of $C_i$.

\mn
We now show that $\Gamma$ does not contain any Grassmannian subspaces $S$ of diameter $3$ or more. We do this by creating a certain set of points $\cC_S$ inside $S$ and showing that it must be contained in a set $\cC$ of points in $\Gamma$, which is in fact smaller, thus obtaining a contradiction.
\ble\label{lem:no A5,3}
The geometry $\Gamma$ does not contain subspaces of type $A_{m,l}$ for
 $3\le l\le m-2$.
\ele
\pf
If such a subspace exists, then there also exists a subspace $S$ of type $A_{5,3}$.
Let $x,y\in S$ with $d_S(x,y)=3$.
We consider the sets
$$\begin{array}{rl}
\cC_S&=\{u\in S\mid d_S(x,u)=2, \mbox{ and }d_S(u,y)=1\},\\
\cC&=\{u\in S\mid \dist_\Gamma(x,u)=2g,2q, \mbox{ and }\dist_\Gamma(u,y)=1\}.\\
\end{array}$$
By Lemma~\ref{lem:A3,2 in Dn,n-2}, for any $u\in \cC_S$ we have $\dist_\Gamma(u,x)=2g,2q$, so that $\cC_S\sbe \cC$. Our aim is to derive a contradiction from the fact that $\cC$ is ``too small'' to contain $\cC_S$.

Note that, by Lemma~\ref{lem:subspace distances}, $2\le \dist_\Gamma(x,y)\le 3$.

First let $\dist_\Gamma(x,y)=3$, then by Lemma~\ref{lem:Dk+2}, $x$ and $y$ are contained in a subspace $R$ of type $D_{5,3}$.
Moreover, every $u\in\cC$ is on a geodesic from $x$ to $y$, and since $R$ is convex by Lemma
~\ref{lem:residues are convex}, we may from now on assume that $n=5$ and $k=3$.

Checking the $D_5$ apartment we find that $x$ and $y$ can be in one of three possible relative positions
 (see Figure~\ref{fig:distances in D53}).
Two of these cases (3q and 3hh) enforce $\cC=\emptyset$ contradicting $\cC_S\ne\emptyset$.
In the third case (3h), all points of $\cC$ lie on a line, whereas $\cC_S$ contains triples of pairwise non-collinear points, again a contradiction.

Next, assume $\dist_\Gamma(x,y)=2$.

First assume $\dist(x,y)=2s$.
We know that if $u\in \cC$, then $\dist(x,u)=2g$ or $2q$.
Suppose that $u$ has distance $2g$ to $x$.
Then, $y\cap x^\perp\sbe u\sbe x^\perp\cap y^\perp$. This leaves exactly two possibilities for $u$.
Now suppose that $u$ has distance $2q$ to $x$. Then,
 $u=\langle x\cap y^\perp, t\rangle_W$, where $t$ is some $1$-space in $y-x^\perp$ (since $\dist(x,u)\ne 2s$).
All points but one on the line given by $x\cap y^\perp$ and $\langle x\cap y^\perp,y\rangle_W$ satisfy this condition.
The one remaining point is the unique common neighbor of $x$ and $y$, which is not contained in $S$.
So $\cC$ is the union of two isolated points and a line minus a point. Since $\cC_S$ contains a plane, we cannot have $\cC_S\sbe\cC$. This rules out the case $\dist(x,y)=2s$.

Next assume that $\dist(x,y)=2q$ and that $u\in \cC$.
Again, $\dist(x,u)=2g$ or $2q$.
If $\dist(x,u)=2g$, then $u\sbe x^\perp\cap y^\perp$. In particular, $u\cap y\sbe x^\perp\cap y=y\cap x$, which
 rules out $\dist(x,u)=2g$.
Thus $\dist(x,u)=2q$.
This means that $u=\langle x\cap y, t\rangle_W$, where $t$ is some $1$-space in $y^\perp-x^\perp$.
Since $y^\perp\cap x^\perp$ forms a grid,  if $u'$ is another such $1$-space, then
  $u$ and $u'$ are either not collinear, or if they are, then the line $uu'$ contains a common neighbor of
   $x$ and $y$.
In particular, $\cC$ does not contain any projective planes, but $\cC_S$ does.
This rules out the case $\dist(x,y)=2q$ altogether.

Finally assume $\dist(x,y)=2g$. Let $u\in \cC$. Then, $u\cap y$ contains $x\cap y$ or not.
If not, then $\dist(x,u)\ne 2$ as we can see in Figure~\ref{fig:distances in D53}, so
 $x\cap y\sbe u\cap y$.
If $u\cap x$ is an $(n-3)$-object, then $u$ is a common neighbor of $x$ and $y$, a contradiction.
So $u\cap y=x\cap y$. This leaves two cases: $\dist(x,u)\in \{2s,2g\}$, but as $u\in \cC$, this means
  $\dist(x,u)=2g$.
In particular, $x\cap y\sbe u\sbe x^\perp\cap y^\perp$. In the residue of $x\cap y$ we see that $\cC$ is 
  the subset of a geometry of type $A_{3,2}$ of $2$-spaces missing the $2$-space corresponding to $x$
   and disjoint from the $2$-space corresponding to $y$.
Thus every line meeting $\cC$ in at least two points, in fact has a point outside $\cC$.
Since $\cC_S$ contains full lines only, we cannot have $\cC_S\sbe \cC$.
\qed

\bth\label{thm:parabolic}
A subspace of $\Gamma$ isomorphic to $A_{m,l}$ with $2\le \min\{l,m+1-l\}$ is parabolic.
\eth
\pf
By Lemma~\ref{lem:no A5,3} we must have $l\in\{2,m-1\}$. Up to a graph automorphism we may assume
 $S\cong A_{m,2}$ where $m\ge 3$.
If $m=3$ the claim is Lemma~\ref{lem:A3,2 in Dn,n-2}.
Now assume $m\ge 4$.
Consider points $x, y, z\in S$ such that $y$ and $z$ are collinear, and $\dist(x,w)=2$ for all $w$ on the line spanned by $y$ and $z$.
In $S$ we easily see that the collection of common neighbors to $x$, $y$ and $z$ forms a line.

Now we consider $x$, $y$, and $z$ as points of $\Gamma$ and show that
 they are incident to a common object of type $n$ or $(n-1)$.
By Lemma~\ref{lem:A3,2 in Dn,n-2}, $x$ and $y$ are in a subspace of type $A_{3,2}$, that is parabolic. That is, they are incident to a flag of type $\{n-3\}$, $\{n-4,n\}$, or
 $\{n-4,n-1\}$.
The configuration of the points $x$, $y$, and $z$ cannot be realized in a single $A_{3,2}$-subspace. 
Therefore $x$, $y$, and $z$ cannot all share the same $(n-3)$-object.
However, if $x$ and $y$ do share an $(n-3)$-object that is not also shared with $z$, then Figure~\ref{fig:distances in D53} shows that $x$ has distance $2q$ to $y$, but distance $2s$ to $z$.
Our configuration precludes this.
Therefore we may assume without loss of generality that $x$ and $y$ are incident to a necessarily unique object $u$ of type $n$, but not to a common object of type $n-1$ or $n-3$.
Interchanging the roles of $y$ and $z$, we find that $z$ must share an object $u'$ of type $n$ or $n-1$ with $x$.
We find that $u=\langle x,y\rangle_W$ and $u'=\langle x,z\rangle_W$ are  orthogonal maximal totally singular subspaces in $W$.
It follows that $u=u'$.

Now consider a graph $(\Theta,\sim)$, where $\Theta=\{\{x,y\}\mid x,y\in S, d_S(x,y)=2\}$ in which $\{x,y\}\sim\{x',y'\}$ if, possibly after switching the roles of $x$ and $y$, we have
 $x=x'$ and $y$ is collinear to $y'$ and $x$ is not collinear to any point on the line $yy'$.
We claim that the graph $(\Theta,\sim)$ is connected.
First note that if $\{x,y\}, \{x,y'\}\in \Theta$, then either $x$ has distance $2$ to all points on the line $yy'$, or, since $\dim(V)=m+1\ge 5$, there exists some point $y''$ meeting $\langle x,y,y'\rangle_V$ only in the $1$-space $y\cap y'$.
Now $\{x,y\}\sim\{x,y''\}\sim\{x,y'\}$.
Thus, in order to show connectedness, we can remove the condition on $\sim$ that $x$ not be collinear to any point on the line $yy'$.
Evidently the resulting graph, and hence $(\Theta,\sim)$ itself, is connected.
Moreover, every point $x$ is clearly part of some pair in $\Theta$.
Therefore, by the preceding, all points $x$ lie on the same object of type $n$ (or $n-1$).
\qed
\section{An application}
Theorem~\ref{mainthm1} has several applications. One was suggested by M.~Pankov.
We'll illustrate this with the following observation.

Let $\De$ be a building of type $M_n(\FF)\ne D_n(\FF)$ and let $\Gamma$ be its $k$-Grassmannian as in Table~\ref{table:polar types}.
Let $\Pi$ and $\DP$ be the polar and dual polar space associated to $\De$.
\medskip

\samepage{We first list a result which is fairly well-known.
\ble\label{lem:dual polar autos}
\
\begin{itemize}
\AI{a} $\Aut(\Pi)\cong \Aut(\De)$.
\AI{b} $\Aut(\DP)\cong \Aut(\De)$.
\end{itemize}
\ele
}\pf
(a) Since every $i$-object of $\De$ can be identified as a singular subspace of dimension $i$
 in $\Pi$~\cite{BuSh1974,Ti1986}, we have a homomorphism $\varphi\colon \Aut(\Pi)\to\Aut(\De)$.
Also, for every shadow space $\Theta$ of $\De$ we have
 a homomorphism $\theta\colon \Aut(\De)\to\Aut(\Theta)$.
This holds in particular for every polar Grassmannian $\Theta$ and even more specifically for $\Theta=\Pi,\Gamma,\DP$.
Taking $\Theta=\Pi$ we see that $\theta=\varphi^{-1}$.

(b)
By part (a) we have a homomorphism $\theta\colon \Aut(\De)\to\Aut(\DP)$.
As seen in e.g.~\cite{Ca1982}, the $i$-objects $X$ of $\De$ correspond bijectively to convex closures of pairs of points in $\DP$ at distance $n-i$ from each other in the collinearity graph of $\DP$. This allows to construct the building $\De$ from $\DP$.
Since automorphisms preserve distances, subspaces and convexity, we have a homomorphism $\varphi\colon \Aut(\DP)\to\Aut(\De)$ such that $\varphi\after\theta=\id$.
\qed

\bco\label{cor:autos}
If $M_n\ne D_n, B_{4,2},\wD_{5,2}$, then $\Aut(\Gamma)\cong\Aut(\De)$.
\eco
\pf
For $k=1,n$, this is Lemma~\ref{lem:dual polar autos}.
So we may assume that $2\le k\le n-1$.
By Theorem~\ref{mainthm1} the conditions on $\Gamma$ imply that for every maximal singular subspace $B\le W$ of $\De$ the subspace $S=S(B,0)$ is parabolic and any subspace isomorphic to $S$ is necessarily of the form $S(B',0)$ for some maximal singular subspace $B'\le W$ of $\De$.
Thus there is a bijection between the point set of $\DP$ and the collection of subspaces
 of type $A_{n-1,k}$ in $\Gamma$.

Next, we consider collinearity.
Let $B_1,B_2$ be distinct points of $\Theta$.
Note that $B_1\cap B_2$ is a singular $m$-space in $W$ for some $m=0,1,2,\ldots,n-1$.
For $m\ge k$, this is visible in $\Gamma$ from the fact that $S(B_1\cap B_1,0)=S(B_1,0)\cap S(B_2,0)$ is a non-empty Grassmannian
 subspace of type $A_{m-1,k}$.
In fact, whenever $m\ge k$, the subspace $S(B,0)$ contains $S(B_1\cap B_2,0)$ if and only if $B\spe B_1\cap B_2$.
In particular the subspaces $L$ of type $A_{n-2,k}$ of the form $S(B_1,0)\cap S(B_2,0)$ correspond bijectively to the lines of $\DP$.

Any automorphism of $\Gamma$ preserves the collection of subspaces $S$ of type $A_{n-1,k}$ and also preserves the collection
 of intersections $L$ of type $A_{n-2,k}$ of a pair of such $A_{n-1,k}$-type subspaces.
Therefore we have a homomorphism $\varphi\colon\Aut(\Gamma)\to\Aut(\DP)$ and the result now follows from Lemma~\ref{lem:dual polar autos}.
\qed

\bibliographystyle{alpha}

\begin{flushleft}
Rieuwert J. Blok \\ Department of Mathematics and Statistics\\
Bowling Green State University \\ Bowling Green, OH 43403 \\[2pt]
Tel.: +1 (419) 372-7455 \\
FAX : +1 (419) 372-6092 \\
Email: blokr@member.ams.org\\
\end{flushleft}

\begin{flushleft}
Bruce N. Cooperstein\\
Mathematics Department\\
University of California at Santa Cruz\\
Santa Cruz, CA 95064\\[2pt]
Tel.: +1 (831) 459-2150\\
FAX : +1 (831) 459-3260\\
Email: coop@ucsc.edu
\end{flushleft}
\end{document}